\documentclass[12pt,a4paper]{article}

\usepackage[greek,american]{babel}
\usepackage{amssymb}
\usepackage{amsfonts}
\usepackage{subfigure}
\usepackage[dvips]{graphicx}

\begin{document}
\begin{center}

{{{\Large\sf\bf Distributional Results Relating to the Posterior of a Dirichlet Process Prior}}}\\

\vspace{0.5cm}
{\large\sf Spyridon J. Hatjispyros $^*$, Theodoros Nicoleris$^{**}$ and Stephen G. Walker$^{***}$}\\

\end{center}

\centerline{\sf $^{*}$ Department of Mathematics, University of the Aegean,}
\centerline {\sf Karlovassi, Samos, GR-832 00, Greece.} 
\centerline{\sf $^{**}$ Department of Economics, National and Kapodistrian University of Athens,}
\centerline{\sf Athens, GR-105 59, Greece. }  
\centerline{\sf $^{***}$Department of Mathematics, University of Texas at Austin,}
\centerline{\sf Austin, Texas 7812, USA. }

\begin{abstract} The aim of this paper is to find distributional results for the posterior parameters which arise in the Sethuraman (1994) representation of the Dirichlet process. These results can then be used to derive simply the posterior of the Dirichlet process.

\vspace{0.1in} \noindent {\sl Keywords:} Bayesian nonparametrics; Dirichlet process; Distributional theory.
\end{abstract}

\vspace{0.5in} \noindent {\bf 1. Introduction.} The aim of this paper is to provide a self contained exposition of the posterior of a Dirichlet process prior. In particular we concentrate on distributional results relating to the posterior of the parameters which arise in the Sethuraman (1994) representation of the Dirichlet process. This then provides a self contained proof of the posterior of the Dirichlet process also being a Dirichlet process.

While the derivation of the posterior Dirichlet process has been done before, notably in Ferguson (1973), we provide the derivation based on the constructive definition of the Dirichlet process given in Sethuraman (1994). This representation is in the form
$P=\sum_{j=1}^\infty w_j\delta_{\theta_j}$, where the $\theta=(\theta_j)$ are independent and identically
distributed (i.i.d.) from some distribution $G$ and the $(w_j)$ are the weights. Specifically,
$w_1=v_1$, and for $j>1$, $w_j=v_j\prod_{l<j}(1-v_l)$ and the $v=(v_j)$ are i.i.d. beta$(1,c)$ for some $c>0$. We find the posterior for $(\theta,v)$ given a single observation from $P$, say $X=x$, and then find distributional results for these variables which provide the posterior for the Dirichlet process.  

The posterior distribution based on this representation was also covered by Sethuraman. His target was the posterior distribution of $P|[X=x]$, as is ours, yet we work with the posterior distribution of $[v,\theta|X=x]$. With this we can work directly with Bayes theorem. With key distributional results relating to this posterior, which are of interest in their own right, we are able to derive the posterior for $P$. The results are all self contained and a clear exposition of the posterior emerges. 

\vspace{0.2in} \noindent {\bf 2. The posterior for $[\theta,v|X=x]$.} To start, if $X$ is a random sample from $P$ then it must take the value of one of the $(\theta_j)$, and takes the value $\theta_j$ with probability $w_j$. Hence, it is natural to introduce the latent variable $J$ to indicate which component $X$ came from. Hence, we have that  
$$
\mbox{P}(J=j|\,\theta,v)=w_j,\quad \mbox{P}(x|J=j, \theta,v)=1(x=\theta_j).$$
If we use $\pi(\theta)$ and $\pi(v)$ to denote the prior distributions for $\theta$ and $v$, respectively, then Bayes Theorem gives
\begin{equation}
\label{rel1}
\pi(\theta,v|x)=\sum_{j=1}^\infty{\pi(\theta,v,x,J=j)\over h(x)}
=\sum_{j=1}^\infty{\pi(\theta)\pi(v)w_j1(x=\theta_j)\over h(x)},
\end{equation}
where $h(x)$ is the density function corresponding to the distribution $G$. 
Now since
$$
{\pi(\theta)1(x=\theta_j)\over h(x)}=
\prod_{l\ge 1}h(\theta_l){1(x=\theta_j)\over h(\theta_j)}=
\prod_{l\neq j}h(\theta_l)1(x=\theta_j),
$$
equation (\ref{rel1}) becomes 
\begin{equation}
\label{rel2}
\pi(\theta,v|x)=\sum_{j=1}^\infty {w_j\pi(v)\over \bar{w}_j}
\prod_{l\neq j} h(\theta_l) 1(x=\theta_j)\, \bar{w}_j,
\end{equation}
where $\bar{w}_j=\mbox{E}(w_j)$, so
$$\bar{w}_j={1\over c+1}\left({c\over c+1}\right)^{j-1}$$
for $j\geq 1$.

Hence, from equation (2), it is clearly seen that $\pi(\theta,v|x)$ is a mixture
and that with probability $\bar{w}_j$, $\theta$ is distributed as 
$${\bf 1}(\theta_j=x)\,\prod_{l\ne j}h(\theta_l)$$
and $v$ is distributed as 
$$\bar{w}_j^{-1} v_j\prod_{l<j}(1-v_l)\,\pi(v).$$
Therefore, given $J$, the $(\theta_j)$ are independent, as are the $(v_j)$, and for the latter,
for each $j$, define the density $\pi_j(v_1^j,v_2^j,\ldots)$ by
$$v_l^j \sim \left\{ \begin{array}{ll}
\mbox{beta}(1,c+1) & l<j \\ \\
\mbox{beta}(2,c) & l=j \\ \\
\mbox{beta}(1,c) & l>j.\end{array}\right.$$
So the distribution of $v$ given $j$ is $\pi_j(v_1^j,v_2^j,\ldots)$ which is clearly seen by considering
$\pi(v|j)\propto w_j\,\pi(v)$. 
Also easy to see is $\pi(\theta|j,x)$, so define $\pi_j(\theta_1^j,\theta_2^j,\ldots)$ by $\theta_l^j\sim G$ for $l\ne j$ and $\theta_j^j=x$.

This completes the posterior for $[\theta,v|x]$. 
To sample from the posterior of $[P|x]$; sample $j$ with probability $\mbox{P}(j)=\bar{w}_j$. Then take
$$P_j=\sum_{l\ne j} w_l^j \delta_{\theta_l^j} + w_j^j \delta_x,$$
which we can write as
\begin{equation}
\label{Pj}
P_j=(1-w_x) \sum_{l=1}^\infty q_l^j \delta_{\theta_l} + w_x\delta_x,
\end{equation}
where $w_x=w_j^j$, the $(\theta_l)$ here are i.i.d. $G$, and
$$q_l^j=\left\{\begin{array}{ll}
w_l^j/(1-w_x)  & l<j \\ \\
w_{l+1}^j /(1-w_x) & l\geq j.\end{array}\right.$$
Here, we have put $w_x=w_j^j$ and we will also define the random weights $q$ via: $q=(q_l^j)_{l=1}^\infty$ with probability $\bar{w}_j$.

In order to get a posterior representation we need to show a few distributional results based on the posterior of $[\theta,v|x]$. The important results are as follows: We want to show that, marginally, i.e. integrating out $J$,  
\begin{description}
\item 1. $w_x$ is beta$(1,c)$.

\item 2. $q$ is a set of Dirichlet process weights with parameters $(c,G)$. 

\item 3. $w_x$ and $q$ are independent.

\end{description}

\noindent We will go through these results in order. First, it is easy to see that $w_x$ is beta$(1,c)$.

\vspace{0.1in} \noindent {\bf Theorem $1$.} For all $n\geq 1$ it is that $\mbox{E}(w_x^n)=\mbox{E}(v^n)$, where $v$ is beta$(1,c)$.

\vspace{0.1in} \noindent {\sc Proof.} Now
$$\begin{array}{ll}
\mbox{E}(w_x^n) & = \sum_{j=1}^\infty \bar{w}_j \frac{\Gamma(2+c)\Gamma(n+2)}{\Gamma(n+2+c)}\left(\frac{c+1}{c+1+n}\right)^{j-1} \\ \\
& = \sum_{j=1}^\infty \frac{\Gamma(2+c)\Gamma(n+2)}{\Gamma(n+2+c)}\frac{1}{c+1}\left(\frac{c}{c+1+n}\right)^{j-1} \\ \\
& = \frac{\Gamma(1+c)\Gamma(n+2)}{\Gamma(n+2+c)}\frac{c+1+n}{n+1} \\ \\
& =\mbox{E}(v^n) , \end{array}$$
thus completing the proof.

\vspace{0.2in} \noindent
To establish the result for $q$, we first define
$$h_l^j=q_l^j/(1-q_1^j-\cdots-q_{l-1}^j)$$
and put $h_l=h_l^j$ with probability $\bar{w}_j$.
Clearly,
$$h_l^j=\left\{\begin{array}{ll}
v_{l}^j/[1-\psi_{j,l}(1-v_l^j)] & l<j \\ \\
v_{l+1}^j  & l\geq j, \end{array}\right.$$
where $\psi_{j,l}=v_j^j\prod_{l<r<j}(1-v_r^j)$.
We now want to show that the $(h_l)$ are mutually independent and each $h_l$ is beta$(1,c)$. This will establish that $q$ is a set of random Dirihclet process weights.  

\vspace{0.1in} \noindent {\bf Theorem $2$.} Marginally $h_l$ is beta$(1,c)$.

\vspace{0.1in} \noindent {\sc Proof.} Now
$$\mbox{P}(h_l\geq y)   =\sum_{j=1}^\infty \mbox{P}(h_l^j\geq y)\,\bar{w}_j$$
which can be written as
$$(1-y^c)\sum_{j=1}^l \frac{1}{1+c}\left(\frac{c}{1+c}\right)^{j-1}+\sum_{j>l}\mbox{P}\left(\frac{v_l^j}{1-\psi_{j,l}(1-v_l^j)}\geq y \right) \frac{1}{1+c}\left(\frac{c}{1+c}\right)^{j-1}.$$
The first term is easily seen to be
$$(1-y)^c \left(1-(c/(1+c))^l\right)$$
and the second term is given by
$$\begin{array}{l}
 \sum_{j>l} \mbox{E}\left\{ \frac{(1-y)^{c+1}}{(1-\psi_{j,l}y)^{c+1}}\right\}\frac{1}{1+c}\left(\frac{c}{1+c}\right)^{j-1} \\ \\
= (1-y)^{c+1} \sum_{j>l} \frac{1}{1+c}\left(\frac{c}{1+c}\right)^{j-1}\sum_{n=0}^\infty \small {\left(\begin{array}{l} c+1 \\ n \end{array}\right)}y^n \mbox{E}\left(\psi_{j,l}^n\right).\end{array}$$
But
$$\mbox{E} \left(\psi_{j,l}^n\right)=\frac{\Gamma(2+c)\Gamma(n+2)}{\Gamma(n+2+c)}\left(\frac{c+1}{c+1+n}\right)^{j-l-1}$$
so the second term becomes, after some few lines of algebra,
$$(1-y)^c \left(c/(1+c)\right)^l.$$
Hence the result when putting first and second terms together.

\vspace{0.1in} \noindent We now look at independence. We will first need the following lemma

\vspace{0.1in} \noindent {\bf Lemma $1$.} The distribution of $[\psi_{J,l}|J>l]$ is ${\rm beta}(1,c)$.

\vspace{0.1in} \noindent {\sc Proof.} We have that 
$$
\mbox{E}\left(\psi_{J,l}|J>l\right)\,=\,
{1\over \mbox{P}\{J>l\}}\mbox{E}\left(\psi_{J,l}{\bf 1}_{\{J>l\}}\right)\,=\,
\left({c\over c+1}\right)^{-l}\sum_{j=l+1}^\infty\bar{w}_j\mbox{E}(\psi_{J,l}{\bf 1}_{\{J>l\}}|J=j).
$$
Using the fact that for $j>l$ we have
\begin{eqnarray}
\nonumber
\mbox{E}(\psi_{J,l}{\bf 1}_{\{J>l\}}|J=j) & = &
\mbox{E}((v_j^j)^n|j>l)\, \mbox{E}((1-v_{l+1}^j)^n|j>l)\cdots \mbox{E}((1-v_{j-1}^j)^n|j>l)\\
\nonumber
 & = & {\Gamma(n+2)\Gamma(c+2)\over\Gamma(c+n+2)}\left({c+1\over c+n+1}\right)^{j-1-l}
\end{eqnarray}
then
\begin{eqnarray}
\nonumber
\mbox{E}(\psi_{J,l}|J>l) & = &
{1\over c+1}\left({c\over c+1}\right)^{-l}{\Gamma(n+2)\Gamma(c+2)\over\Gamma(c+n+2)}
\sum_{j=l+1}^\infty {\left(c\over c+1\right)}^{j-1}{\left(c+1\over c+n+1\right)}^{j-l-1}\\
\nonumber
 & = & {\Gamma(n+1)\Gamma(c+1)\over\Gamma(c+n+1)}\,=\,\mbox{E}(v^n),
\end{eqnarray}
where $v\sim{\rm beta}(1,c)$, which proves lemma $1$.

We now put in a useful and necessary Lemma.

\vspace{0.1in} \noindent {\bf Lemma $2$.} If $v\sim{\rm beta}(1,c+1)$ and $\psi\sim{\rm beta}(1,c)$
are independent random variables then 
\begin{enumerate}
\item $\xi=\psi (1-v)\sim{\rm beta}(1,c+1)$. 
\item $\zeta=v/(1-\psi(1-v))\sim {\rm beta}(1,c)$.
\item The random variables $\xi$ and $\zeta$ are independent.
\end{enumerate}

\vspace{0.1in} \noindent {\sc Proof.} Straightforward.

\vspace{0.1in}
Before showing the mutual independence of the $(h_l)$, we first establish the independence of $h_l$ and $w_x$ for any $l$. This can be helpful as a first step and is also the last of our required distributional results as well. 

\vspace{0.1in} \noindent {\bf Theorem $3$.} It is that, for every $l$ and $n,m\geq 1$,
$$\mbox{E}(h_l^n\,w_x^m)=\mbox{E}(h_l^n)\,\mbox{E}(w_x^m).$$

\vspace{0.1in} \noindent {\sc Proof.} Now
$$
\mbox{E}(h_l^n\,w_x^m)  =\sum_{j=1}^l \mbox{E}((v_{l+1}^j)^n) \mbox{E}((w_j^j)^m)\bar{w}_j + \mbox{E}(\zeta^n)\,\mbox{E}\left\{\xi^m\prod_{r<l}(1-v_r)^m\right\}\bar{w}_{j>l}
$$
where $v$ is beta$(1,c+1)$, $\psi$ is beta$(1,c)$, $(v_r)$ are i.i.d. beta$(1,c+1)$
and $\bar{w}_{j>l}=\sum_{j>l}\bar{w}_j$. 
For $j\leq l$ it is that $v_{l+1}^j$ is beta$(1,c)$ and it is also easy to show that
$v/(1-\psi(1-v))$ is also beta$(1,c)$. Therefore,
$$\mbox{E}(h_l^n\,w_x^m) =\mbox{E}(v^n) \times k(m)$$
where $v$ is now beta$(1,c)$ and $k$ some function of $m$. Hence, $h_l$ and $w_x$ are marginally independent.

\vspace{0.1in} \noindent Now we need to show that $h_l$ and $h_{l'}$ are independent. The idea for proving this is essentially the same as for showing that $h_l$ and $w_x$ are independent. Let us assume without loss of generality that $l>l'$. 

\vspace{0.1in} \noindent {\bf Theorem $4$.} It is that $h_l$ and $h_{l'}$ are independent.

\vspace{0.1in} \noindent {\sc Proof.}
Now, for $j=1,\ldots,l$, it is that $h_l^j=v_{l+1}^j$ and each $v_{l+1}^j$ is beta$(1,c)$ and is independent of $h_{l'}^j$. Now, for $j>l$,
$$h_l^j=\frac{v_l^j}{1-\psi_{j,l}(1-v_l^j)}$$
and
$$h_{l'}^j=\frac{v_{l'}^j}{1-\psi_{j,l}(1-v_l^j)\prod_{l'\leq r<l}(1-v_r^j)}. $$
Given $j>l$, and based on marginals derived from previous results, we have
$$h_l|[j>l]=\frac{v}{1-\psi(1-v)}$$
and
$$h_{l'}|[j>l]=\frac{w}{1-\psi(1-v)\phi }$$
where $v$ and $w$ are independent beta$(1,c+1)$, and independent of $\psi$ which is beta$(1,c)$, and $\phi$ is independent of both $\psi$ and $v$. So $h_l|[j>l]$ and $h_{l'}|[j>l]$ are independent, due to Lemma 2, and $h_l|[j>l]$ is beta$(1,c)$. Combining this with the results for $h_l[j\leq l]$ which is beta$(1,c)$ and also independent of $h_{l'}|[j\leq l]$, we have the required result.

\vspace{0.2in} \noindent {\bf 3. Posterior for $[P|X=x]$.} With the distributional results for $[\theta,v|X=x]$ we can now easily find the posterior for $[P|X=x]$.  

\vspace{0.2in} \noindent {\bf Lemma $3$.} Relation (\ref{Pj}) marginally becomes 
$P_x=(1-\beta)P+\beta\delta_x$ where $\beta\sim\mbox{beta}(1,c)$, $P\sim\mbox{DP}(c,G)$
and $P_x\sim\mbox{DP}(c+1,{cG+\delta_x\over c+1})$

\vspace{0.1in} \noindent {\sc Proof.} The fact that relation (\ref{Pj}) marginally is
$P_x=(1-\beta)P+\beta\delta_x$ comes from theorems $1$, $2$, and, $3$.
Now to show that $P_x$ is a sample from the Dirichlet process with concentration
parameter $c+1$ and base measure $c\,G+\delta_x\over c+1$,
we consider a finite partition of the support $X$ of $G$,
$X=A_1\cup\cdots\cup A_k$ and that $x\in A_j$ for some $1\le j\le k$ then
$$
(P_x(A_1),\ldots,P_x(A_k))\stackrel{d}{=}
\beta(\delta_x(A_1),\ldots,\delta_x(A_k))+(1-\beta)(P(A_1),\ldots,P(A_k)).
$$
The random vector $(\delta_x(A_1),\ldots,\delta_x(A_k))$ has the degenerate Dirichlet 
distribution $D(e_{j1},\ldots,e_{jk})$ with $e_{ji}=1(j=i)$. On the other hand
$$
(P(A_1),\ldots,P(A_k))\stackrel{d}{=}D(c\,G(A_1),\ldots,c\,G(A_k))
$$
and
$$
\mbox{beta}(1,c)=\mbox{beta}\left(\sum_{i=1}^k e_{ji},\sum_{i=1}^k c\,G(A_i)\right),
$$
then from standard theory we have
$$
(P_x(A_1),\ldots,P_x(A_k))\stackrel{d}{=}
D(c\,G(A_1)+\delta_x(A_1),\ldots,c\,G(A_k)+\delta_x(A_k))
$$
which gives the desired result.

\vspace{0.2in} \noindent {\bf 4. Discussion.} Our style of proof for the posterior Dirichlet process obviously most closely resembles the one provided in Sethuraman (1994). However, there are notable differences. These are the key distributional results for $[\theta,v|X=x]$ which we believe are necessary for the rigorous derivation of the posterior and which are of interest in their own right. From these results we can also see why the Dirichlet process is unique. We have also provided a framework for exploring alternative stick--breaking priors, whereby the $(v_l)$ have alternative beta distributions. 

\vspace{0.2in} \noindent {\bf References.} 

\begin{description}

\item Ferguson, T.S. (1973). A Bayesian Analysis of some Nonparametric Problems.
{\sl Annals of Statistics} {\bf 1}, 209–-230.

\item Sethuraman, J. (1994). A Constructive Definition of Dirichlet Priors. 
{\sl Statistica Sinica} {\bf 4}, 639--650.

\end{description}

\end{document}